\documentclass[10pt]{amsart}
\usepackage[dvipdfm]{graphicx,color}
\usepackage{url}

\numberwithin{equation}{section}
\newtheorem{theorem}{Theorem}[section]
\newtheorem{lemma}[theorem]{Lemma}
\theoremstyle{definition}

\newtheorem{corollary}[theorem]{Corollary}

\theoremstyle{remark}

\numberwithin{equation}{section}

\def \bfo {\begin {eqnarray*} }
\def \efo {\end {eqnarray*} }
\def \ba {\begin {eqnarray*} }
\def \ea {\end {eqnarray*} }
\def \beq {\begin {eqnarray}}
\def \eeq {\end {eqnarray}}

\def \bfo {\begin {displaymath} }
\def \efo {\end {displaymath} }
\def \beq {\begin {eqnarray}}
\def \eeq {\end {eqnarray}}
\def \ba {\begin {eqnarray*}}
\def \ea {\end {eqnarray*}}

\def\tilde{\widetilde}

\def \bfo {\begin {displaymath} }
\def \efo {\end {displaymath} }
\def \beq {\begin {eqnarray}}
\def \eeq {\end {eqnarray}}
\def \ba {\begin {eqnarray*}}
\def \ea {\end  {eqnarray*}}

\def \R {{\mathbb {R}}}

\def \supp {\hbox{supp }}

\def \p {\partial}

\def \tilde{\widetilde}
\def \proofbox {$\square$\medskip}

\def\supp{{\mathop{\rm{supp}}}}

\begin{document}

\title[Recent progress of inverse scattering theory]{Recent progress of inverse scattering theory\\ on non-compact manifolds}

\author[Hiroshi Isozaki]{Hiroshi Isozaki}

\address{Hiroshi Isozaki, Institute of Mathematics \\
University of Tsukuba,
Tsukuba, 305-8571, Japan;
Yaroslav Kurylev, Department of Mathematics \\
University College of London, United Kingdom;
Matti Lassas, Department of Mathematics and Statistics\\
University of Helsinki, Finland
}
\email{isozakih@math.tsukuba.ac.jp,y.kurylev@ucl.ac.uk,mjl@rni.helsinki.fi}

\author[Yaroslav kurylev]{Yaroslav Kurylev}
\address{}
\email{}

\author[Matti Lassas]{Matti Lassas}
\address{}
\email{}

\thanks{}

\maketitle
\begin{center}
{\small{\it Dedicated to Gunther Uhlmann on the occasion of his 60th birthday}}
\end{center}

\begin{abstract}
We give a brief survey for the recent development of  inverse scattering theory on non-compact Riemannian manifolds. The main theme is the reconstruction of the manifold and the metric from the scattering matrix. 

\end{abstract}


\section{Introduction}
Scattering theory  for waves in  classical or quantum physics has a long history, and nowadays 
there is an abundance of works devoted to the forward and inverse problems of potential 
scattering for Schr{\"o}dinger equations, and obstacle scattering  for  wave equations. 
Much less is known, however, about the inverse scattering on manifolds, where the main theme 
is the reconstruction of the manifold itself and its Riemannian metric from the scattering 
matrix associated with the Laplace-Beltrami operator. We have been working on this subject by choosing typical examples possessing  characteristic features of the problem. 
The aim of this 
paper is to present a birds-eye view based on the results we have obtained so far.

In \S2, we discuss a basic framework of our issue. Scattering of waves is a time-dependent 
phenomenon, however, it is a consequence of properties  of the continuous spectrum of the 
underlying Laplacian. Therefore, we formulate the problem mainly from the stationary view point 
with emphasis on the solution space of the {\it Helmholtz equation}.
In \S 3, we explain the outline of our theory for the forward and inverse problems ignoring  
the detailed assumptions of the metric. 
In \S4, we pick up four metrics we have already studied. They are well-known standard examples, 
and, viewed from the growth order of volume near infinity, range from an exponentially growing 
case to an exponentially decaying case. 
In \S5, we devote ourselves to the case of 2-dim. arithmetic surface and discuss the inverse scattering from cusp.

In \S6, we mention some recent results for metrics with intermediate behavior.

Due to the lack of space, we cannot talk about all of the important works. For example, \cite
{BorJuPer05}, \cite{Borth}, \cite{Tanya09}, \cite{MU}, \cite{PSU},
\cite{Perry07}  are dealing with problems related to ours from different view points.


\section{Review of forward and inverse problems}

\subsection{Time-dependent scattering problem}
Let us consider ${\mathbb R}^n, \, n>1,$ with a Riemannian metric $G = g_{ij}dx^idx^j$,
where we use the Einstein summation convention of omitting sum signs. We assume that
asymptotically, when $|x| \to \infty$, $g_{ij}(x) \simeq \delta_{ij},$, see e.g. (\ref{S4EuclideAssump}).
How can one recover the metric tensor $g_{ij}$ from some "physical" observations? 
One way is to consider  wave motions associated with this metric:
\begin{equation}
\partial_{t}^2v(t,x) =\Delta_Gv(t,x),
\label{S1Waveeq}
\end{equation}
where $\Delta_G$ is the (negative) Laplace-Beltrami operator associated with $G$.
In the remost past (i.e. $t \to -\infty$) and in the remote future (i.e. $t\to\infty)$, the solution to  (\ref{S1Waveeq}) behaves like $v_{\pm}$ satisfying $\partial_t^2v_{\pm}= \Delta v_{\pm}$, $\Delta$ being the standard Euclidean Laplacian on ${\mathbb R}^n$.  The mapping 
assigning the far-field pattern of $v_-$ to that of $v_+$ is called the {\it scattering operator}. It is a common belief that, under natural situations,  this scattering operator determines the original physical system, i.e.\ the metric, at least up to changes of coordinates.
In other words, from the observation of all  far-field patterns of waves at infinity, one can reconstruct the metric $G$. 
For elliptic inverse boundary value problems on compact or non-compact Riemannian manifolds there are also non-uniqueness results; one can use  the fact that the
measurements are invariant in diffeomorphisms of the manifold  to create counterexamples for the uniqueness of inverse problems and even invisibility results,
see \cite{GLU1,GKLU0,Ke,LTU} and the references in
 \cite{GKLU1,GKLU2}

Another way is to use the wave motion of quantum mechanical particles, in which case we use the Schr{\"o}dinger equation
\begin{equation}
i \partial _t v(t,x) = \Delta_Gv(t,x),
\label{S1TimedepSch}
\end{equation}
instead of (\ref{S1Waveeq}).   (See e.g. \cite{ReSi} or \cite{Ya}).
Both of these time-dependent scattering problems are reformulated in the same stationary (time-independent) picture. 

\subsection{Stationary scattering theory} Consider a time-periodic motion of (\ref{S1Waveeq}), i.e. put $v(t,x) = e^{-ikt}u(x)$ in (\ref{S1Waveeq}). Then we get
\begin{equation}
(- \Delta_G - k^2)u(x)=0, \quad k > 0.
\label{S1StationarEq}
\end{equation}
Since all solution of the Helmholtz equation (\ref{S1StationarEq}) with unperturbed 
Laplacian, $\Delta$,
can be written as  superpositions of plane waves and taking into account that
$g_{ij} \simeq \delta_{ij}$,
  we seek a solution of (\ref{S1StationarEq}) 
admitting the following asymptotic expansion
\begin{equation}
u (x)\simeq e^{ik\omega\cdot x} + \frac{e^{ikr}}{r^{(n-1)/2}}a(\hat x,\omega;k), \quad \hat x = x/r, \quad 
{\rm as} \quad r = |x| \to \infty,
\label{S1StationarAsymp}
\end{equation}
where $\omega \in S^{n-1}$. The 2nd term of the right-hand side represents the scattered spherical wave. In the case of quantum mechanics, for a given beam of particles with initial direction $\omega$, $|a(\theta,\omega;k)|^2$ is proportional to the number of particles scattered to the $\theta$-direction. 
This is  the physical quantity observed in the experiment. 
Let $A(k):L^2(S^{n-1})\to L^2(S^{n-1})$ be the integral operator with the kernel 
$ a(\hat x,\omega;k)$. 
By a suitable choice of the constant $C(k)$, $S(k) = I -C(k) A(k)$ is a unitary operator 
on $L^2(S^{n-1})$, which is the well-known {\it Heisenberg's S-matrix}. The above mentioned scattering operator in the time-dependent formulation is written in terms this S-matrix.

\subsection{Geometric scattering theory} 
Let us look at  the stationary scattering theory from a geometrical view point.
 Consider the restriction of the Fourier transform on $S^{n-1}$:
\begin{equation}
\left(\mathcal F_0(k)^{\ast} \varphi\right)(x) = 
\int_{S^{n-1}}e^{ik\omega\cdot x}\varphi(\omega)d\omega, \quad k >0.
\nonumber
\end{equation}
This is sometimes called the {\it Herglotz integral} or the {\it Poisson integral}. It is the eigenoperator of $- \Delta$ in the sense that it satisfies the equation
$$
(- \Delta - k^2)\mathcal F_0(k)^{\ast}\varphi = 0.
$$
By the stationary phase method, it admits the asymptotic expansion:
\begin{equation}
(\mathcal F_0(k)^{\ast}\varphi)(x) \simeq C_+(k)\frac{e^{ikr}}{r^{(n-1)/2}}\varphi(\hat x) + 
 C_-(k)\frac{e^{-ikr}}{r^{(n-1)/2}}\varphi(-\hat x) , \quad r=|x|\to\infty.
\label{S1F0kAsympto}
\end{equation}
Let us note that the spherical waves $e^{\pm ikr}/r^{(n-1)/2}$   
 appear in the asymptotic expansion of the Green functions  $ G_0^{(\pm)}(x,y;k)$ 
of $- \Delta -k^2$ in ${\mathbb R}^n$:
$$
\int_{{\mathbb R}^n} G_0^{(\pm)}(x,y;k)f(y)dy \simeq 
\widetilde C_{\pm}(k)\frac{e^{\pm ikr}}{r^{(n-1)/2}}\widehat f(\pm k\hat x),
$$
where $\widehat f(\xi)$ is the Fourier transform of $f(y)$.
The unit sphere $S^{n-1}$ can be regarded as a boundary at infinity of ${\mathbb R}^n$. Therefore, for ${\mathbb R}^n$ with Euclidean metric, one can associate the {\it manifold at infinity}, $S^{n-1}$, and the integral transform on it, $\mathcal F_0(k)^{\ast}$. It is an eigenoperator of $- \Delta$ and has the asymptotic expansion (\ref{S1F0kAsympto}).

The stationary scattering theory asserts that these properties are transferred to ${\mathbb R}^n$ with metric $G$. Namely, one can associate a generalized eigenoperator $\mathcal F(k)^{\ast}$ on $L^2(S^{n-1})$ satisfying
$$
(-\Delta_G - k^2)\mathcal F(k)^{\ast}\varphi = 0, \quad {\rm for} \ {\rm all}\  k^2 \in \sigma_c(-\Delta_G),
$$
and admitting the asymptotic expansion (\ref{S1F0kAsympto}) with $\mathcal F_0(k)^{\ast}$ replaced by $\mathcal F(k)^{\ast}$. Moreover, by introducing a Banach space $\mathcal B^{\ast}$ by
\begin{equation}
\mathcal B^{\ast} \ni u \Longleftrightarrow 
\sup_{R>1}\frac{1}{R}\int_{|x|<R}|u(x)|^2dx < \infty,
\label{S2BastRn}
\end{equation}
one can show that
\begin{equation}
\mathcal N(k) := \{u \in \mathcal B^{\ast}\, ; (- \Delta_G - k^2)u = 0\} = \mathcal F(k)^{\ast}\left(L^2(S^{n-1})\right),
\nonumber
\end{equation}
and that for any $\varphi_- \in L^2(S^{n-1})$, there exist unique $u \in \mathcal N(k)$ and $\varphi_+ \in L^2(S^{n-1})$ satisfying
\begin{equation}
u (x)\simeq C_+(k)\frac{e^{ikr}}{r^{(n-1)/2}}\varphi_+(\hat x) + 
 C_-(k)\frac{e^{-ikr}}{r^{(n-1)/2}}\varphi_-(\hat x) , \quad r=|x|\to\infty.
\nonumber
\end{equation}
Here, by $f\simeq g$, we mean
$$
\lim_{R\to\infty}\frac{1}{R}\int_{|x|<R}|f(x)-g(x)|^2dx=0.
$$
The mapping
\begin{equation}
S_{geo}(k) : L^2(S^{n-1}) \ni \varphi_- \to \varphi_+ \in L^2(S^{n-1})
\nonumber
\end{equation}
is unitary.
In fact, it is related with the above S-matrix as follows:
\begin{equation}
S_{geo}(k) = JS(k),
\label{SgeoJ}
\end{equation}
where $J\phi(\omega) = \phi(-\omega)$.

Motivated by the above, let us give an overview how a geometric scattering theory
on a Riemannian manifold can be formulated.
We are given a non-compact Riemannian manifold $\mathcal M$ having the boundary at infinity $M$, which is a compact Riemannian manifold of dimension $n-1$. We then construct an operator $\mathcal F(k)^{\ast} : L^2(M) \to L^2_{loc}(\mathcal M)$ such that
$$
(- \Delta_G - k^2)\mathcal F(k)^{\ast}\varphi = 0, \quad {\rm for} \ {\rm all} \ k^2 \in \sigma_c(-\Delta_G),
$$
with the property that, by a suitable choice of the Banach space $\mathcal B^{\ast}$ of functions
in $\mathcal M$,
$$
\mathcal N(k) := \{ u \in \mathcal B^{\ast}\, ; \, (-\Delta_G-k^2)u=0\} = 
\mathcal F(k)^{\ast}\left(L^2(M)\right).
$$
Moreover, there exist functions $w_{\pm}(k) \in \mathcal B^{\ast}$ such that for any $\varphi_- \in L^2(M)$ there exist unique $u \in \mathcal N(k)$ and $\varphi_+ \in L^2(M)$ admitting the asymptotic expansion
$$
u \simeq w_+(k)\varphi_+ + w_-(k)\varphi_-, 
$$
and the mapping
$$
S_{geo}(k) : L^2(M) \ni \varphi_- \to \varphi_+ \in L^2(M)
$$
is unitary. We then formulate our inverse problem as follows :
{\it
Reconstruct the manifold $\mathcal M$ and its metric from the knowledge of the S-matrix $S_{geo}(k)$.}

\subsection{Inverse scattering}
The inverse potential scattering on the half-line was solved in 1950's by Gel'fand-Levitan and Marchenko \cite{GeLe}, \cite{Mar55}. The multi-dimensional extension was proposed by the seminal work of Faddeev \cite{Fa76}. We do not enter into the detailed exposition here (see e.g. \cite{Iso03}), however, let us emphasize  that the so called {\it CGO solution} of the Schr{\"o}dinger equation introduced by Calder{\'o}n and implemented into a powerfull tool by Sylvester and Uhlmann \cite{SylUhl87} (in fact, it was done independently of Faddeev) stimulated the works of Nachman \cite{Na} and Khenkin-Novikov \cite{KheNov} to bring the progress in the multi-dimensional inverse problem.

As for our inverse scattering problem, we use the {\it boundary control} (BC) method.
It is now regarded as an established theory for  the inverse spectral problem, however,  a long time was necessary to settle the whole idea.
The precursor of the BC method appeared in the work of M. G. Krein on the 1-dimensional wave equation \cite{Kr51a}, \cite{Kr51b}. Although there is a similarity of Krein's theory to that of Gel'fand-Levitan-Marchenko, his idea is based on the hyperbolic nature of the wave equation, i.e. finite velocity of the wave propagation and
the notion of the domain of influence. However, by passing to the Fourier transform in $t$, 
he used it in the form of the analyticity properties with respect to the spectral parameter. Blagovestcenskii \cite{Bla71a} analyzed Krein's idea, using the finite velocity and  controllability of the filled domain, to derive a Volterra-type equation for unknown functions. In \cite{Bla71b}, he also found a crucial idea of evaluating the inner product of waves with given data on the boundary by the spectral data. These ideas were then extended to the multi-dimensional case by the works of Belishev and  Belishev-Blagovestcenskii \cite{Be87}, \cite{BeBla92}. Belishev and Kurylev \cite{BeKu92} used the BC method to solve the Gel'fand inverse problem on compact Riemannian manifolds. We must also mention that the step of controllability was completed by Tataru \cite{Tata}, which gave a final form for Holmgren's uniqueness theorem.
The BC method can also be used to reconstruct non-smooth manifolds \cite{AKKLT}, \cite{KLY} or solve
inverse problems for the heat and Schr\"odinger equations on Riemannian manifolds \cite{KKLM}.
For  the more detailed history of the BC method and the related reconstruction methods, see \cite{Be97,Bingham,KruL,KKL01}.

\section{Bird's-eye view of the inverse scattering problem for the metric}
We shall discuss here "typical" results which we expect to hold. Therefore, we state Theorems and Lemmas "formally" without specfying the detailed assumptions.

\subsection{General framework}
We consider an $n$-dimensional, non-compact, connected, Riemannian manifold $\mathcal M$ with the following properties :
\begin{equation}
\mathcal M = \mathcal K \cup \mathcal M_1\cup\cdots\cup\mathcal M_N,
\label{S3M=KcupMi}
\end{equation} 
where $\mathcal K$ is relatively compact, and $\mathcal M_i$'s (called {\it ends}) are non-compact. We allow $\mathcal K$ to have an arbitrarily metric, however, each $\mathcal M_i$ is assumed  to be diffeomorphic to $(1,\infty)\times M_i$, where $M_i$ is a compact $(n-1)$-dimensional manifold endowed with a metric $h_{i}=(h_i)_{\alpha\beta}(y)dy^\alpha dy^\beta$.  Let us denote points
of $(1,\infty)\times M_i$ by $(r,y)$.
We also assume that on each end $\mathcal M_i$, the metric $G$ is asymptotically equal to
the warped product:
\begin{equation}
G \sim (dr)^2 + \rho_i(r)(h_i)_{\alpha\beta}(y)dy^\alpha dy^\beta, \quad {\rm as} \quad r \to \infty,
\label{S3Warpedprod}
\end{equation}
with some positive function $\rho_i(r) \in C^{\infty}({\mathbb R})$. We assume $\mathcal M$ to be complete, hence $- \Delta_G\big|_{C_0^{\infty}(\mathcal M)}$ is essentially self-adjoint. For the sake of simplicity, assume that $\sigma_c(-\Delta_G) = [E_0,\infty)$, and put
$$
H = -\Delta_G - E_0.
$$
The properties of the continuous spectrum of $H$ are determined by those of ends $\mathcal M_i$. Therefore, we reduce our analysis to a suitable model space, whose choice depends largely on the nature of $\rho_i(r)$, and shall be discussed separately.

\subsection{Analysis on the model space}
 To fix our assumptions,
we  consider the case when the boundary $\partial \mathcal M$ is empty. 
Note that,
however, one could consider, in a similar manner, also the case when 
$\partial \mathcal M \neq \emptyset$ 
with Neumann, Dirichlet or other  self-adjoint  boundary conditions.
Given $\rho(r) \in C^{\infty}({\mathbb R})$, we take 
$\mathcal M_0 = I_0\times M_0$, where $I_0 = (c_0,\infty)$ is an open interval such 
that $c_0 < 1$, equipped with the metric $G_0 = (dr)^2 + \rho(r)h_0$. 
 By imposing the Neumann boundary condition at $r=c_0$, we assume that 
the Laplace-Beltrami operator $H_0 = - \Delta_{G_0}$  on $\mathcal M_0$ is self-adjoint in 
$L^2(I_0,L^2(M_0);\rho(r)^{(n-1)/2}dr)$. 
Assume that $\sigma_{c}(-\Delta_{G_0}) = [E_0,\infty)$, and put 
$H_0 = - \Delta_{G_0}-E_0$,  $R_0(z) = (H_0-z)^{-1}$. Let us call $H_0$ 
the {\it free operator} on $\mathcal M_0$. Our first task is to prove the  
limiting absorption principle, which is going to be explained below.

Let $A = \int_{\mathbb R}\lambda dE_A(\lambda)$ be a self-adjoint operator on a Hilbert space $\mathcal H$. Then for $\lambda \in \sigma(A)$, $(A-\lambda)^{-1}$ does not exist. However, for $\lambda \in \sigma_c(A)$, sometimes one can prove the existence of the limit
$$
(A- \lambda \mp i0)^{-1} = \lim_{\epsilon\to0}(A - \lambda \mp i\epsilon)^{-1} 
\in 
{\bf B}(\mathcal H_+ ;\mathcal H_-),
$$
where ${\bf B}(X;Y)$ means the set of all bounded operators from a Banach space $X$ to a Banach space $Y$, and $\mathcal H_{\pm}$ are suitable Banach spaces skirting $\mathcal H$:
$$
\mathcal H_+ \subset \mathcal H \subset \mathcal H_-,
$$
with inclusions dense and continuous, and $\mathcal H_- = (\mathcal H_+)^{\ast}$ through the inner product of $\mathcal H$. This is called the {\it limiting absorption principle}, abbreviated as LAP. If the LAP is proved on an interval $I \subset \sigma_c(A)$, by Stone's formula, we have
$$
E_A(I) \mathcal H \subset \mathcal H_{ac}(A),
$$
where the absolutely continuous subspace $\mathcal H_{ac}(A)$ is defined by
$$
\mathcal H_{ac}(A) \ni u \Longleftrightarrow (E_A(\lambda)u,u) \ {\rm is} \ {\rm absolutely} \ {\rm continuous} \ {\rm w.r.t.} \ d\lambda.
$$
In all examples we have already treated, the space $\mathcal H_-$ is defined as follows.
\begin{equation}
\mathcal H_- \ni u \Longleftrightarrow
\sup_{R>1}\frac{1}{R}\int_{c_0<r<R}\|u(r)\|^2_{L^2(M_0)}\rho(r)^{(n-1)/2}dr < \infty.
\end{equation}
This is a dual space of some Banach space $\mathcal B$ using a dyadic decomposition of the manifold $\mathcal M_0$, and is denoted by $\mathcal B^{\ast}$ below.
To prove LAP, usually we need to avoid the end points of $\sigma_c(H)$, the eigenvalues embedded 
in the continuous spectrum, and thresholds (the energies at which the nature of scattering changes). 
Let us call the set of these points the exceptional set, and denote it by $\mathcal E_0$. This is 
shown to be discrete with possible accumulation points $0$ and on thresholds.

Our next task is to analyse the asymptotic behavior of the resolvent $R_0(\lambda \pm i0)$ as $r \to \infty$. 
 For $f,g\in \mathcal H_-$ we use the notation $f \simeq g$ to denote
$$
f\simeq g \Longleftrightarrow \lim_{R\to\infty}\frac{1}{R}\int_{c_0<r<R}
\|f(r)-g(r)\|^2_{L^2(M_0)}\rho(r)^{(n-1)/2}dr=0.
$$
Then for any $k > 0$ such that $k^2 \not\in \mathcal E_0$, there exist bounded operators 
$\mathcal F_0^{(\pm)}(k) \in {\bf B}(\mathcal B\, ; \, L^2(M_0))$ and 
functions $w_{\pm}(k) \in \mathcal  B^{\ast}$ such that for $f \in \mathcal B$,
\begin{equation}
R_0(k^2 \pm i0)f \simeq w_{\pm}(k)\mathcal F_0^{(\pm)}(k)f.
\label{R0Asympto}
\end{equation}

\subsection{Gluing the estimates for ends}
Having established the LAP and the asymptotic expansion (\ref{R0Asympto}) on each end, we need to glue them together to obtain global results.
Let $\{\chi_j\}_{j=0}^N$ be a partition of unity on $\mathcal M$ such that
$\sum_{j=0}^N\chi_j = 1$,  and, for $j=1,\cdots,N$,  ${\rm supp}(\chi_j) \subset \mathcal M_j$, $\chi_j(r) = 1$ if $r > 2$. We take $\widetilde \chi_j \in C^{\infty}(\mathcal M_j)$  such 
that ${\rm supp}\,(\widetilde \chi_j) \subset \mathcal M_j$ and $\widetilde \chi_j=1$
on ${\rm supp}\,\chi_j$. Let $\mathcal H_{free(j)}$ be the free operator on $\mathcal M_j$
 and assume, for simplicity, that for all $j=1, \dots, N$, we have 
$\sigma_c(\mathcal H_{free(j)})=[E_0, \infty)$. Then 
we put $R_{free(j)}(z) = (H_{free(j)}-z)^{-1}$ and 
$$
A(z) = \sum_{j=1}^NA_j(z)\widetilde \chi_j,
$$
$$
A_j(z) = [H,\chi_j]R_{free(j)}(z) + \chi_j(H-H_{free(j)})\widetilde \chi_jR_{free(j)}(z).
$$
Then we have
\begin{equation}
R(z) = \sum_{j=1}^N\chi_jR_{free(j)}(z)\widetilde\chi_j + R(z)(\chi_0 - A(z)).
\label{S3PerturbationR(z)}
\end{equation}
This formula and the perturbation technique enable us to prove LAP for $H$. 
Let $\mathcal E = \big( \cup_{j=1}^N\mathcal E_j\big)\cup\big(\sigma_p(H)\cap[0,\infty)\big)$, 
$\mathcal E_j$ being the exceptional set for 
$H_{free(j)}$. One can show that $\mathcal E$ is discrete with possible accumulation points 0 and on thresholds of $H_{free(j)}$.
We define the space $\mathcal B^{\ast}$ by gluing ${\mathcal B_j}^{\ast}$ using the partition 
of unity,  namely, 
$$\mathcal B^{\ast}=\{u \in L^2_{loc}(\mathcal M):\, \chi_j u \in \mathcal B^{\ast}_j   \}.
$$
\begin{theorem}
 For $\lambda \in [0,\infty)\setminus\mathcal E$, there exists a limit $R(\lambda \pm i0) \in {\bf B}(\mathcal B\, ;\, \mathcal B^{\ast})$.
\end{theorem}
 
By taking the adjoint in (\ref{S3PerturbationR(z)}), we can derive the asymptotic expansion of the resolvent on each end:

\begin{equation}
R(k^2\pm i0)f \simeq w_j^{(\pm)}(k)\mathcal F_j^{(\pm)}(k)f, \quad f \in \mathcal B.
\label{S3RzAsympto}
\end{equation}
We now introduce the $L^2$-space over the manifold at infinity of $\mathcal M$:
\begin{equation}
{\bf h}_{\infty} = \oplus_{j=1}^N L^2(M_j).
\label{S3hinfty}
\end{equation}
We also put
\begin{equation}
\mathcal F^{(\pm)}(k) = (\mathcal F_1^{(\pm)}(k),\cdots,\mathcal F_N^{(\pm)}k)) 
\in {\bf B}(\mathcal B\,;\,{\bf h}_{\infty}),
\label{S3F(k)define}
\end{equation}
\begin{equation}
\widehat{\mathcal H} = L^2((0,\infty),{\bf h}_{\infty};\,dk).
\label{S3widehatH}
\end{equation}
Here by $L^2(I,{\bf h};dk)$ with given auxiliary Hilbert space $\bf h$, we mean the 
Hilbert space of all ${\bf h}$-valued $L^2$-functions on an interval $I$ with respect to 
the measure $dk$.

\begin{theorem}\label{Thm 3.2}
(1) The operator $\mathcal F^{(\pm)}$ defined for $f \in \mathcal B$ by
$(\mathcal F^{(\pm)}f)(k) = \mathcal F^{(\pm)}(k)f$ is uniquely extended to a partial isometry on $L^2(\mathcal M)$ with initial set $\mathcal H_{ac}(H)$ and final set $\widehat{\mathcal H} $. \\
\noindent
(2) For $f \in D(H)$, we have
$$
(\mathcal F^{(\pm)}Hf)(k) = k^2(\mathcal F^{(\pm)}f)(k).
$$
(3) $\mathcal F^{(\pm)}(k)^{\ast}$ is an eigenoperator in the sense that
$$
(H - k^2)\mathcal F^{(\pm)}(k)^{\ast}\varphi = 0, \quad \hbox{for all }\varphi \in {\bf h}_{\infty}.
$$
(4) For any $f \in \mathcal H_{ac}(H)$, the inversion formula holds:
$$
f = \int_{0}^{\infty}\mathcal F^{(\pm)}(k)^{\ast}(\mathcal F^{(\pm)}f)(k)dk.
$$
\end{theorem}

We have thus constructed a spectral representation (generalized Fourier transformation) for $H$.

\begin{lemma}
For any $k^2 \in (0,\infty)\setminus \mathcal E$,
$$
\mathcal N(k) := \{u \in \mathcal B^{\ast}\, ; \, (- \Delta_G-E_0-k^2)u=0\}
= \mathcal F^{(\pm)}(k)^{\ast}{\bf h}_{\infty}.
$$
\end{lemma}

\begin{theorem} \label{Thm 3.4}
For any $\varphi^{(-)}= (\varphi_1^{(-)},\cdots,\varphi_N^{(-)}) \in {\bf h}_{\infty}$, there exist unique $u \in \mathcal N(k)$ and $\varphi^{(+)}= (\varphi_1^{(+)},\cdots,\varphi_N^{(+)}) \in {\bf h}_{\infty}$ such that on each end $\mathcal M_j$, the asymptotic expansion
$$
u \simeq w_j^{(-)}(k)\varphi_j^{(-)} + w_j^{(+)}(k)\varphi_j^{(+)}
$$
holds.  Moreover the operator
\ba
& &S_{geo}(k) : {\bf h}_{\infty} \to {\bf h}_{\infty},\\
& &S_{geo}(k) \varphi^{(-)} = \varphi^{(+)}
\ea
is unitary.
\end{theorem}

\subsection{Rellich's theorem, Radiation condition}
In the Euclidean case, that is, for the manifold $\mathbb R^n$, there are two types of boundary 
conditions at infinity. A solution to the equation $(- \Delta_G - k^2)u = f, \, k>0,$  in 
${\mathbb R}^n$  is said to satisfy the outgoing radiation condition, if $(\frac{\partial}{\partial r}-ik)u = o(|x|^{-(n-1)/2})$. Similarly, $u$ is said 
to satisfy the incoming radiation condition, 
if $(\frac{\partial}{\partial r}+ik)u= o(|x|^{-(n-1)/2})$. 
Using the space ${\mathcal B}^{\ast}$, we can rewrite these conditions as
\begin{equation}
\label{radiation}
\lim_{R\to\infty}\frac{1}{R}\int_{1<r<R}\Big\|(\frac{\partial}{\partial r} \mp ik)
u(r)\Big\|_{L^2(S^{n-1})}^2 r^{(n-1)}dr = 0
\end{equation}
Moreover, the classical Rellich theorem says that if $u \in \mathcal B^{\ast}(\Bbb R^n)$
solves a homogeneous Helmholtz equation and
satisfies an incoming or outgoing radiation condition (\ref{radiation}), then $u=0$.

These results can be  extended to a more general case. Namely,  for any $k > 0$ such that $k^2 \not\in\mathcal E$, it is possible to define an incoming and outgoing radiation conditions for
 $u \in \mathcal B^{\ast}$. They are defined by
\begin{equation}
\lim_{R\to\infty}\frac{1}{R}
\sum_{j=1}^N \left(\int_{1<r<R}\Big\|(\frac{\partial}{\partial r} - iE_j(k))
u(\cdot, r)\Big\|_{L^2(M_j)}^2\rho_j(r)^{(n-1)/2}dr \right) = 0,
\label{S3Outgoingrad}
\end{equation}
and
\begin{equation}
\lim_{R\to\infty}\frac{1}{R}
\sum_{j=1}^N \left(\int_{1<r<R}\Big\|(\frac{\partial}{\partial r} + iE_j(k))
u(\cdot, r)\Big\|_{L^2(M_j)}^2\rho_j(r)^{(n-1)/2}dr \right) = 0,
\label{S3Incomingrad}
\end{equation}
corrrespondingly. Here $E_j(k)$ is some operator in $L^2(M_j)$.

In typical situations, the space $\mathcal B^{\ast}$ enjoys the following property.

\begin{theorem}
Let $u \in {\mathcal B}^{\ast}$ satisfies
$(- \Delta_G - E_0 - \lambda)u = 0$ for $\lambda \in (0,\infty)\setminus\mathcal E$
together with an outgoing (\ref{S3Outgoingrad}), or incoming (\ref{S3Incomingrad})
radiation condition.
Then $u = 0$.
\end{theorem}

Namely, as for the decay at infinity, $\mathcal B^{\ast}$ is the smallest space for non-trivial solutions to the Helmholtz equation $(-\Delta_G-E_0-\lambda)u=0$. This is a generalization of the classical Rellich's theorem.

\begin{theorem}
A solution $u \in \mathcal B^{\ast}$ of the equation $(- \Delta_G-E_0 - k^2)u = f\in \mathcal B$ satisfying either the outgoing radiation condition,
or the incoming radiation condition 
is unique and is given by $u = R(k^2 + i0)f$, or $u = R(k^2 -i0)f$, respectively.
\end{theorem}

\subsection{Inverse problem}
Having completed the forward problem, we can now enter into the inverse problem. 
Since $\mathcal M$ has $N$ number of ends, the S-matrix $S_{geo}(k)$ is an $N\times N$ matrix-valued unitary operator, $S_{geo}(k) = \big(S_{ij}(k)\big)_{i,j=1}^N$. We pick up the entry $S_{11}(k)$, and try to reconstruct the whole manifold $\mathcal M$ and the Riemannain metric from $S_{11}(k)$.

Take the end $\mathcal M_1$ and split $\mathcal M$ into 2 parts:
\begin{equation}
\mathcal M = \mathcal M_{ext} \cup \mathcal M_{int}, \quad 
\mathcal M_{ext} = \mathcal M_1 \cap \{r > 2\}, \quad \mathcal M_{int} = \mathcal M\setminus\mathcal M_{ext}.
\label{S3SplitM}
\end{equation}
  Let $H_{ext}$ be the  self-adjoint unbounded  operator $- \Delta_G$ in $L^2(\mathcal M_{ext})$
 identifed by the Neumann boundary condition
at
$r=2$.
Let  $H_{int}$ be the similarly defined  self-adjoint unbounded  operator $- \Delta_G$ in
 $L^2(\mathcal M_{int})$
with Neumann boundary condition at
$r=2$.
 One can then solve the boundary value problem on $\mathcal M_{ext}$ and $\mathcal M_{int}$ :
\begin{equation}
\left\{
\begin{split}
& (-\Delta_G-E_0-k^2)u = 0, \quad {\rm in} \quad \mathcal M_{ext}, \\
& \frac{\partial u}{\partial r}= f, \quad {\rm on} \quad  \partial\mathcal M_1=\{r=2\},\\
& u \ {\rm satisfies} \ {\rm the} \ {\rm outgoing} \ {\rm radiation} \ {\rm condition}.
\end{split}
\right.
\label{S3ExteriorProblem}
\end{equation}
\begin{equation}
\left\{
\begin{split}
& (-\Delta_G-E_0-k^2)u = 0, \quad {\rm in} \quad \mathcal M_{int}, \\
& \frac{\partial u}{\partial r}= f, \quad {\rm on} \quad  \partial\mathcal M_1=\{r=2\},\\
& u \ {\rm satisfies} \ {\rm the} \ {\rm outgoing} \ {\rm radiation} \ {\rm condition} \ {\rm if} \ {\rm N \geq 2}.
\end{split}
\right.
\label{S3InteriorProblem}
\end{equation}
Let $\Lambda_{ext}(k)$ and $\Lambda_{int}(k)$ be the associated Neumann to Dirichlet maps:
\begin{equation}
\Lambda_{ext}(k) : f \to u_{ext}\big|_{\{r=2\}},
\label{S3ExtNDmap}
\end{equation}
\begin{equation}
\Lambda_{int}(k) : f \to u_{int}\big|_{\{r=2\}},
\label{S3IntNDmap}
\end{equation}
$u_{ext}$ and $u_{int}$ being a unique solution to (\ref{S3ExteriorProblem}) and (\ref
{S3InteriorProblem}), respectively. 
 Here and below, we assume that $k^2$ is not on the exceptional sets associated
with $\mathcal H_{ext}$ and $\mathcal H_{int}$ or an eigenvalue of $H_{int}$ if $\mathcal H_{int}$ is a bounded region.

The following lemma is a bridge between the scattering problem and the {\it interior} boundary value problem. 


\begin{lemma}\label{SmatrixandDNmap}
Given $\Lambda_{ext}(k)$ for all $k > 0$, the S-matrix $S_{11}(k)$ and the ND map $\Lambda_{int}(k)$ determine each other for all $k > 0, k^2\not\in\mathcal E$.
\end{lemma}

This is a well-known fact. The proof for the case of the Euclidean space is given in \cite{IsaNa93}, and the hyperbolic space case is deal with in \cite{Iso04}. In the proof of this lemma, Rellich's uniqueness theorem plays a crucial role.

 With the aid of Lemma \ref{SmatrixandDNmap}, assuming that we know the end $\mathcal M_1$, we can obtain $\Lambda_{int}(k)$ for all $k \in {\mathbb C}$ from the $(1,1)$-component of the S-matrix $S_{11}(k)$ for $k > 0$. 
 One can then apply the BC method to $\mathcal M_{int}$, namely we have

\begin{theorem}\label{BCMainTh}
Suppose we are given two (possibly non-compact) Riemannian manifolds $\mathcal M^{(1)}_{int}$, $\mathcal M^{(2)}_{int}$ such that there is a diffeomorphism $\psi:\p  \mathcal M^{(1)}_{int}\to \p \mathcal M^{(2)}_{int}$
and the Dirichlet-to-Neumann maps satisfy $\psi^*(\Lambda_{int}^{(2)}(k)f)=\Lambda_{int}^{(1)}(k)\psi^*f$
for all $k\in (a,b)$, $b>a>0$ and $f\in C^\infty_0(\p  \mathcal M^{(2)}_{int}).$
Then 
$\mathcal M^{(1)}$ and $\mathcal M^{(2)}$ are isometric.
\end{theorem}

 The detailed proof of this theorem it is too long to review  here but we explain below the sketch of the main steps of the proof. The details are seen in \cite{BeKu92,KKL01} for the case of compact manifolds and in \cite{IK12,LSU} for non-compact manifolds. 
 The first step needed to prove Theorem \ref{BCMainTh} is 
  to 
 consider  the hyperbolic problem corresponding to (\ref{S3InteriorProblem})
 \begin{equation}
\left\{
\begin{split}
& (\p_t^2-\Delta_G)u^F(x,t) = 0, \quad {\rm in} \quad \mathcal M_{int}\times \R_+, \\
& \frac{\partial u^F}{\partial \nu}= F(x,t), \quad {\rm on} \quad  \partial\mathcal M_{int}\times \R_+,\\
&u^F(x,t)|_{t=0}=0,\quad \p_tu^F(x,t)|_{t=0}=0.
\end{split}
\right.
\label{S3InteriorProblem BC}
\end{equation}
Assume next that we are given $\partial\mathcal M_{int}$ and the operator 
$ \Lambda_{int}(k)$ for $k\in (a,b)$. 
Since $\Lambda_{int}(k)$ is meromorphic operator depending  on $k \in {\mathbb C}$,
this data determine $ \Lambda_{int}(k)$ for all $k \in {\mathbb C}$.
Since the Fourier transform  $(\mathcal F_{t\to \omega}{u^F})(x,\omega)$ of $u^F(x,t)$ with respect to $t$ satisfies (\ref{S3InteriorProblem})
with $k^2+E_0=\omega^2$,
we can determine the hyperbolic Dirichlet-to-Neumann map $ \Lambda^{(h)}: F\mapsto 
u^F|_{\partial\mathcal M_{int}\times \R_+}$.
This map  can be used to compute inner products of waves.

\begin {lemma}(Blagovestchenskii identity)
\label {l:4.10}
Let $F,H\in C^\infty_0 (\p \mathcal M_{int}\times \R_+)$ and $T>0$.
Then
\beq
\int_{\mathcal M_{int}} u^F(x,T) {u^H(x,T)}\,dV_G(x)=
\label {4.60}
\eeq
\bfo
=\frac{1}{2} \int_L \int_{\p \mathcal M_{int}} (
F(x,t)  {(\Lambda^{(h)} H)(x,s)}-(\Lambda^{(h)}H)(x,t)  {h(x,s)})\,dS_G(x)dtds,
\efo
where
\bfo
L=\{ (s,t):\ 0\leq t+s\leq 2T,\ t<s ,\ t,s>0 \}.
\efo
\end {lemma}

%
%
%
%

\noindent
{\bf Proof.} Let
\bfo
w(t,s)=\int_ {\mathcal M_{int}} u^F(x,t) {u^H(x,s)}\,dV_G(x).
\efo
Then, by integration by parts, we see that
\ba
(\p^2_t-\p^2_s)w(t,s)=
&=&
\int_{\mathcal M_{int}}  [\Delta_Gu^F(x,t) {u^H(x,s)}- u^F(x,t) 
{\Delta_Gu^H(x,s)}]\,dV_G(x)
\\&=&\int_{\p {\mathcal M_{int}} } [\p_\nu u^F(t) {u^H(s)}- u^F(t) 
{\p_{\nu} u^H(s)}]\,dS_G(x)
\\&=&\int_{\p {\mathcal M_{int}} } [\Lambda^{(h)}F(x,t) {u^H(x,s)}- u^F(x,t) 
{\Lambda^{(h)}H(x,s)}]\,dS_G(x).
\ea
Moreover,
\bfo
\left. w \right| _{t=0}=\left. w\right| _{s=0}=0,\quad
\left. \p _t w\right| _{t=0}= \left. \p_s w\right| _{s=0}=0.
\efo
Thus, $w$ is the solution of the initial-boundary value problem
for the one-dimensional wave equation in the domain $(t,s) \in
[0,2T]\times [0,2T]$ with known source and zero initial and $\Lambda^{(h)}$.
Solving this problem,
we determine $w(t,s)$ in the domain where $t+s\leq 2T$ and $t<s$. In particular,
$w(T,T)$ gives the assertion.
\hfill\proofbox

The other result is based on the following fundamental theorem by 
D. Tataru \cite{Tata,Tata2}.

\begin{theorem}
\label{th:tataru}
Let $u=u(x,t)$ solve the wave equation
$u_{tt}-\Delta_G u=0$  in ${\mathcal M_{int}}\times \R$ and
$u|_{\Gamma \times (0,2T_1)} = \p_{\nu}u|_{\Gamma \times (0,2T_1)} =0$,
where $\emptyset \neq \Gamma \subset \p {\mathcal M_{int}}$ is open. 
Then
\beq
\label{2.9p}
u(x,t)=0 \,\, \hbox{in} \, \, K_{\Gamma,T_1},
\eeq
where  
\ba 
K_{\Gamma,T_1}  =\{(x,t)\in {\mathcal M_{int}}\times \R: d(x,\Gamma) < T_1-|t-T_1|\}
\ea
is the double cone of influence.
\end{theorem}
(The proof of this theorem, in full generality, is in \cite{Tata}.
A simplified proof for the considered case is in \cite{KKL01}.)

The observability Theorem \ref{th:tataru} gives rise to the following approximate controllability:

\begin {corollary}
\label{cor:2.4}
For any open $ \Gamma \subset \p M$ and $T_1>0$,
\bfo
\hbox{cl}_{L^2({\mathcal M_{int}})}\{u^F(\cdot,T_1): F \in C^{\infty}_0(\Gamma \times (0,T_1)) \}
= L^2({\mathcal M_{int}}(\Gamma,T_1)).
\efo
Here 
\bfo
{\mathcal M_{int}}(\Gamma,T_1) = \{x\in {\mathcal M_{int}}:\ d(x,\Gamma) < T_1\} = K_{\Gamma,T_1} \cap \{t=T_1\}
\efo is the domain of influence of $\Gamma$ at time $T_1$
and $ L^2({\mathcal M_{int}}(\Gamma,T_1))=\{a\in L^2({\mathcal M_{int}}):\ \supp(a)\subset
{\mathcal M_{int}}(\Gamma,T_1)\}$.
\end{corollary}

\begin{lemma}\label{Main F2} Let $T>0$
and $\Gamma_j\subset \p {\mathcal M_{int}}$, $j=1,\dots,J$,
 be  non-empty, relatively compact open sets, $0\leq T_j^-<T_j^+\leq T$.
Assume we are given $\p {\mathcal M_{int}}$ and the response
operator $\Lambda^{(h)}$. 
  This
data determines the inner product 
\ba
 J^T_{\mathcal N}(F_1,F_2)=\int_{\mathcal N}u^{F_1}(x,t)u^{F_2}(x,t)\,dV_G(x)
\ea
for given $t>0$ and
$F_1,F_2\in C^\infty_0(\p {\mathcal M_{int}}\times \R_+)$,
 where 
 \beq\label{N set}
 {\mathcal N}=\bigcap_{j=1}^J ({\mathcal M_{int}}(\Gamma_j,T_j^+)\setminus {\mathcal M_{int}}(\Gamma_j,T_j^-))\subset {\mathcal M_{int}}.
 \eeq
\end{lemma}

 We give the proof of Lemma \ref{Main F2}  only   for the case when $J=1$ and $T_1^-=0$ 
so that ${\mathcal N}={\mathcal M_{int}}(\Gamma_1,T_1^+)$. 
We note that when the manifold $\p{\mathcal M_{int}}$ and
 the response
operator $\Lambda^{(h)}$ are given, we can determine 
 the metric $G_{\p {\mathcal M_{int}}}$ of the boundary ${\p {\mathcal M_{int}}}$.
We say that the manifold $\p{\mathcal M_{int}}$ with its metric and the response
operator $\Lambda^{(h)}$ form the {\it boundary data}.

Consider first the case $F_1=F_2=F$. Let $B=\Gamma_1\times [T-T_1^+,T]).$
For all $H\in C^\infty_0(B)$ it holds  that $\supp(u^H(\cdotp,T))\subset {\mathcal N}$,
and thus
\ba
& &\|u^F(T)-u^H(T)\|_{L^2({\mathcal M_{int}})}^2
\\
& &=\int_{\mathcal N} (u^F(x,T)-u^H(x,T))^2dV_G(x)
+\int_{{\mathcal M_{int}}\setminus {\mathcal N}} (u^F(x,T))^2dV_G(x).
\ea
Let $\chi_{\mathcal N}(x)$ be the characteristic function of the set $\mathcal N$.
By Corollary \ref{cor:2.4}, there is $H\in C^\infty_0(B)$ such
that the norm $\|\chi_{\mathcal N}u^F(T)-u^H(T)\|_{L^2(M,dV_\mu)}$
is arbitrarily small. This shows that $ J^T_{\mathcal N}(F_1,F_2)$ can be found
by
\beq\label{eq: minimize 0}
 J^T_{\mathcal N}(F,F)=\|u^F(T)\|_{L^2({\mathcal M_{int}})}^2-
   \inf_{H\in C^\infty_0(B)} \mathcal F(H),
 \eeq  
where
\ba
\mathcal F(H)=\|u^F(T)-u^H(T)\|_{L^2({\mathcal M_{int}})}^2.
\ea
Since $\mathcal F(H)$ can be computed when  the  boundary data is given by Lemma \ref{l:4.10},
it follows that we can determine $ J^T_{\mathcal N}(F,F)$ for any $F\in C^\infty_0(\p M\times \R_+)$. 
Now, since
\ba
 J^T_{\mathcal N}(F_1,F_2)=\frac 14( J^T_{\mathcal N}(F_1+F_2,F_1+F_2)- 
 J^T_{\mathcal N}(F_1-F_2,F_1-F_2)),
\ea
Lemma \ref{Main F2}  follows in the case when $J=1$ and 
$T_1^-=0$.

 To reconstruct $({\mathcal M_{int}},G)$,
we use a special representation, the
{\it boundary distance representation}, $R({\mathcal M_{int}})$, of ${\mathcal M_{int}}$ and 
show that the boundary  data determine $R({\mathcal M_{int}})$.
Consider a map  $R:{\mathcal M_{int}}  \to C(\p {\mathcal M_{int}})$,
\beq
\label{2.1}
R(x)= r_x(\cdot);\quad r_x(z)= d_g(x,z),\,\, z \in \p {\mathcal M_{int}},
\eeq
i.e., $r_x(\cdot)$ is the {\it distance function} from $x\in {\mathcal M_{int}}$ to the points on $\p {\mathcal M_{int}}$. The 
image $R({\mathcal M_{int}}) \subset C(\p {\mathcal M_{int}})$ of $R$ is called the
boundary distance representation of ${\mathcal M_{int}}$. The set $R({\mathcal M_{int}})$ 
is a metric space with the distance inherited from $C(\p {\mathcal M_{int}})$ 
with the standard norm 
$\|f\|_\infty=\max_{z\in \p {\mathcal M_{int}}}|f(z)|$
of $C(\p {\mathcal M_{int}})$. 
We denote this distance by $d_C$. The map $R$, 
due to the triangular inequality, is Lipschitz,
\beq
\label{2.2}
d_C(r_x,r_y) \leq d_g(x,y).
\eeq
Note that if $\p {\mathcal M}_{int}$ is unbounded, which corresponds to the possible case
of unbounded $M_i$, $r_x(\cdot)$ becomes unbounded and, instead of the distance 
$d_C(r_x, r_y)$, it is convenient to use a modified one
\[
d_b(r_x, r_y)= \mathop{\sup_{z\in \p {\mathcal M}_{int}} }
\left(\frac{|r_x(z)- r_y(z)|}{1+|r_x(z)- r_y(z)|}\right).
\]
When $x, y$ run over a compact in ${\mathcal M}_{int}$, $d_b(r_x, r_y)$ is
equivalent to $ d_C(r_x, r_y)$.

The map $R:{\mathcal M_{int}}\to R({\mathcal M_{int}})\subset C(\p {\mathcal M_{int}})$
is an embedding. Many results of differential geometry, such as Whitney or Nash embedding theorems, concern the question
how an abstract manifold can be embedded to some simple space, such as a
higher dimensional Euclidean space. In the inverse problem we need
to construct a "copy" of the unknown manifold in some known space,
and as we assume that the boundary is given, 
we do this by embedding the manifold ${\mathcal M_{int}}$ to the known, although infinite dimensional
function space $C(\p {\mathcal M_{int}})$.

The basic observation needed to construct $R({\mathcal M_{int}})$ is that the set 
$\mathcal N$ given in
(\ref{N set}) has non-zero measure if and only if 
\ba
P(\{(\Gamma_j,T_j^+,T_j^-):\ j=1,\dots, J\})=\sup_f  J^T_N(f,f),
\ea
is non-zero.

\begin{theorem}
\label{th:2.6}
Let $\{z_n\}_{n=1}^{\infty}$ be a dense set on $\p {\mathcal M_{int}}$ and
$r(\cdot) \in C(\p {\mathcal M_{int}})$ be an arbitrary continuous function. Then $r\in R({\mathcal M_{int}})$ if and only if for all $L \in {\Bbb N}_+,$
it holds that
\beq
\label{2.11}
P\Big(\Big\{(z_j,r(z_n)+\frac 1L, r(z_n)-\frac 1L); j=1,\dots,L\Big\}\Big)>0.
\eeq
Moreover, condition (\ref{2.11}) can be verified when
the manifold $\p{\mathcal M_{int}}$ and
 the response
operator $\Lambda^{(h)}$ are given. Hence the pair $(\p{\mathcal M_{int}},\Lambda^{(h)})$
determines uniquely the boundary 
distance representation $R({\mathcal M_{int}})$ of $({\mathcal M_{int}},G)$.
\end{theorem}

The idea of to proof Theorem \ref{th:2.6} is that we see that the set 
$ \{x\in  {\mathcal M_{int}};\  |d(z_j,x)-r(z_n)|<\frac 1L\hbox{ for all }j=1,\dots,L\}$
contains a point $x_0$ for all $L$ if and only if $r_{x_0}(z)=r(z)$ for all $z\in \p {\mathcal M_{int}}$.

\bigskip

Next let  us consider a compact manifold $({\mathcal M_{int}},G)$ that is geodesically
regular manifold, i.e.,
\begin {enumerate}

\item[i)] For any $x,y\in {\mathcal M_{int}}$ there is a unique geodesic $\gamma$
joining these points.
\item[ii)] Any geodesic $\gamma([a,b])$ can be continued to a geodesic
  $\gamma([a',b'])$ whose end-points lie on the boundary.
\end{enumerate}
Consider $R({\mathcal M_{int}})$ as a metric space $(R({\mathcal M_{int}}), d_{\infty})$ with 
the distance
inherited from $L^\infty(\p {\mathcal M_{int}})$,
$
d_\infty(r_x,r_y)=\|r_x-r_y\|_{L^\infty(\p {\mathcal M_{int}})}.
$
Then the mapping $R$ is an isometry, i.e.
\bfo
d_\infty(r_x,r_y)=d(x,y).
\efo
Indeed, let  $\gamma([a,b])$ be the shortest geodesic from $y$ to $x$.
Continue this geodesic  to the shortest geodesic  $\gamma([a',b])$,
where $z=\gamma(a')\in \p {\mathcal M_{int}}.$ Then
\bfo
r_x(z)-r_y(z)=|\gamma([a',b])|-|\gamma([a',a])|=b-a=d(x,y).
\efo
Hence,
\bfo
d(x,y)\leq d_\infty(r_x,r_y).
\efo
The opposite inequality is valid by  the triangular inequality, yielding that
$d(x,y)=d_\infty(r_x,r_y)$.
This means that, if we know {\it a priori} that $({\mathcal M_{int}},G)$
is geodesically regular,
  then $(R({\mathcal M_{int}}), d_{\infty}))$
is isometric to $({\mathcal M_{int}},G)$ and the problem of the reconstruction
of $({\mathcal M_{int}},G)$ from  the data $(\p{\mathcal M_{int}},\Lambda^{(h)})$ is solved. The above ideas can be generalized
to prove Theorem \ref{BCMainTh}.

Let us note that only one diagonal entry  $S_{11}^{(2)}(k)$ is necessary to recover the whole manifold
using Lemma \ref{SmatrixandDNmap} and Theorem \ref{BCMainTh}. 
Therefore, our final inverse procedure is as follows. 

\begin{theorem}\label{S3MainTh}
Suppose we are given two non-compact Riemannian manifolds $\mathcal M^{(1)}$, $\mathcal M^{(2)}$ having $N_1$ and $N_2$ numbers of ends, respectively. Let $\mathcal M_1^{(1)}$ and $\mathcal M_1^{(2)}$ be their first ends, and suppose that there exists $r_0 >1$ such that $\mathcal M_1^{(1)}\cap\{r>r_0\}$ and $\mathcal M_1^{(2)}\cap\{r>r_0\}$ are isometric. Assume furthermore that their $(1,1)$-entries of S-matrices coincide : $S_{11}^{(1)}(k) = S_{11}^{(2)}(k)$ for all $k > 0$, $k^2 \not\in \mathcal E^{(1)}\cup\mathcal E^{(2)}$. 
Then 
$\mathcal M^{(1)}$ and $\mathcal M^{(2)}$ are isometric.
\end{theorem}


\section{Examples}

\subsection{Model metric}
Spectral properties of $- \Delta_G$ depends largely on the growth of the manifold $\mathcal M$ 
at infinity.  Looking at an end $(1, \infty) \times M$,
we pick up here 4 basic examples, and examine what is going on.
 We have so far studied the following cases:

\medskip
\noindent
(1) $G = (dr)^2 + e^{2r}h$, 

\medskip
\noindent
(2) $G = (dr)^2 + r^2h$, 

\medskip
\noindent
(3) $G = (dr)^2 + h$, 

\medskip
\noindent
(4) $G = (dr)^2 + e^{-2r}h$,

\medskip
\noindent
where $h$ is a metric on $M$.

\medskip
\noindent
These are well-known classical examples, and have some distinguished properties. 

\subsection{Tools for the resolvent estimates}
Before going into the details, we explain here the method for resolvent estimates. There are several ways of proving LAP for Laplacians on manifolds. One  is E. Mourre's abstract commutator theory \cite{Mou81}, another is Melrose' theory of scattering metric \cite{MaVa05}, \cite{Me95}. Our method is different from both of them. We employ the classical Eidus' approach of integration by parts \cite{Eid69}.  Given the equation $(-\Delta_{G_0}-E_0-z)u=f$ on the model space $\mathcal M_0 = I_0\times M_0$, we expand $u$ by the eigenvectors of $-\Delta_{M_0}$. The problem is then reduced to the 1-dimensional case. 
The main step is to multiply this 1-dimensional equation by (the derivative of ) solution and integrate by parts to obtain some identities. Standard  arguments from perturbation theory then prove the necessary resolvent estimates. See e.g. \cite{IKL10}, Lemmas 2.4, 2.5, Theorem 2.7, or \cite{IK12}, Lemmas 2.4 $\sim$ 2.8.
This method is no less powerful than the above two machineries, and gives us results by elementary computation.
All of the above 4 examples are treated by this method.

\subsection{Euclidean metric}
Let us start with the Euclidean metric (2). In this case, we take $\mathcal M_0 = {\mathbb R}^n$, and the perturbed metric $g_{ij}dx^idx^j$ is assumed to satisfy
\begin{equation}
\partial_x^{\alpha}(g_{ij}(x) - \delta_{ij}) = O(|x|^{-1-|\alpha|-\epsilon_0}), \quad \hbox{ for all  $\alpha$ and some $\epsilon_0>0$}.
\label{S4EuclideAssump}
\end{equation}

There are  plenty of papers dealing with scattering in the asymptotically Euclidean space, and all the requisites are prepared in e.g. \cite {IK12} or \cite{Ya}. Consulting them, one can see that  the arguments in \S3 work well, and Theorem \ref{S3MainTh} holds.  In this case, there are no embedded eigenvalues in $(0,\infty)$. The radiation condition is the standard one:
$$
\lim_{R\to\infty}\frac{1}{R}\int_{|x|<R}\big|\big(\frac{\partial}{\partial r} \mp ik\big)u\big|^2dx = 0.
$$

\subsection{Cylindrical ends}
Let us consider (3). In the daily life example, it arises as the problem of waveguide. As a model space, 
we take $M_0$ to be a compact $(n-1)$-dimensional manifold without boundary. Let
$\mathcal M_0 = (1,\infty)\times (M_0, h)$, 
and $\Delta_0 =  \Delta_{(M_0,h)}$. The perturbed metric is 
assumed to satisfy
$$
|\partial_X^{\alpha}(g_{ij}(X) - h_{ij}(x))|\leq C_{\alpha}(1 + r)^{-1-\epsilon}, \quad\hbox{for all  $\alpha$ and some $\epsilon>0$}
$$
where $x$ denotes the variable in $M_0$, and $X = (r,x)$.

The waveguide has many features different from the conical Euclidean metric (2). In the first place, it may have eigenvalues embedded in the continuous spectrum. 
The 2nd feature is that the scattering in the waveguide has many channels. Let $\{\lambda_m\}_{m=1}^{\infty}$ be the eigenvalues of $-\Delta_0$. Then, since
$-\Delta_G$ has the form
$$
-\Delta_G = - \big(\frac{\partial}{\partial r}\big)^2 - \Delta_{0},
$$
when the total system has energy $\lambda>0$, only the states with $\lambda_m \leq \lambda$ occur in the scattering phenomena.   The radiation condition is thus defined by
$$
\lim_{R\to\infty} \frac{1}{R} 
\int_{1<r<R}\big\|\big(\frac{\partial}{\partial r} \mp  i\sqrt{\lambda + \Delta_{0}}\big)u(r,\,\cdotp)\|^2_{L^2(M_0)}dr = 0.
$$
Because of this channel property, the scattering matrix changes its size according to the energy, higher the energy, bigger the size of the S-matrix.

To discuss the inverse scattering, we assume that one end is strictly cylindrical. Namely, the end $\mathcal M_1$ is equal to $[1,\infty)\times M_0$ with metric
$(dr)^2 + h$. We take an artificial boundary $r=2$ in $\mathcal M_1$, and split $\mathcal M$ into the exterior and interior domains. As in \S 3, we  derive interior N-D map from the knowledge of $(1,1)$-entry of the S-matrix. 
However, the physical S-matrix is not sufficient to determine the interior N-D map.
The key observation here is that the physical S-matrix $S_{11}(k)$ admits an analytic continuation into the {\it upper} (physical) half plane, and, down to the real axis, 
determines non-physical scattering matrix, which is defined using the exponentially growing solutions of the reduced wave equation instead of the physical plane wave. 
Therefore Theorem \ref{S3MainTh} holds in this case. 

The exponentially growing solution is the crucial idea found by Faddeev in his 
multi-dimensional inverse scattering theory, as well as Calder{\'o}n's inverse boundary value problem. What is interesting is that this apparently artificial exponentially growing solutions appear naturally in the waveguide problem. The details are given in \cite{IKL10}.

We can also allow the ends of type (2) and (3) at the same time. Namely, if our ends $\mathcal M_1,\cdots, \mathcal M_N$ are composed of two parts: $\mathcal M_1, \cdots, \mathcal M_{\mu}$, which are Euclidean ends, and $\mathcal M_{\mu+1},\cdots,\mathcal M_N$, which are cylindrical ends, the results in subsections 4.3, 4.4 are applied as well.


\subsection{Hyperbolic space}\label{S4HypSpace}
We turn to the case (1). The problem here is that in dimensions $\geq 3$, the infinity of hyperbolic manifold may have a complicated structure. Therefore, we restrict ourselves to the simple case, and as a model space, we take $\mathcal M_0 = M_0\times (0,\infty)$ equipped with the metric 
$((dy)^2 + h(x,dx))/y^2$, where $M_0$ is a compact manifold equipped with metric $h(x,dx)$ . This model has two infinities, $y=0$ corresponding to the  infinite volume, which we call {\it regular end}, and $y=\infty$ corresponding to {\it cusp}.  We take $H_0 = -\Delta_{\mathcal M_0} - (n-1)^2/4$. Then $\sigma(H_0) = [0,\infty)$. We split $\mathcal H_0$ into two parts : $\{0<y<1\}$ and $\{y>1\}$. By the change of variable $y = e^{-r}$, the former becomes a model space for (1). The change of variable $y = e^r$ makes the latter to a model space for (4).

The perturbed metric is assumed to have the form
$$
ds^2 = y^{-2}\big((dy)^2 + h(x,dx) + A(x,y,dx,dy)\big),
$$
$$
A(x,y,dx,dy) = \sum_{i,j=1}^{n-1}a_{ij}(x,y)dx^idx^j + 2\sum_{i=1}^{n-1}a_{in}(x,y)dx^idy + a_{nn}(x,y)(dy)^2,
$$
where $a_{ij}(x,y)$ satisfies, for some $\epsilon_0>0$
$$
|\widetilde D_x^{\alpha}D_y^{\beta}\, a(x,y)| \leq C_{\alpha\beta}(1 + |\log y|)^{-\min(|\alpha|+\beta,1)-1-\epsilon_0}, \quad 
{\rm for\ all} \quad  \alpha,\beta.
$$
Here, $D_y = y\partial y$, $\widetilde D_x = \tilde y(y)\partial_x$, $\tilde y(y) \in C^{\infty}((0,\infty))$, $\tilde y(y) = y$ for 
$y>2$, $\tilde y(y) = 1$ for $0 < y < 1$.

If one of the ends is a regular infinity, there are no embedded eigenvalues. If all the ends are cusps, there may be embedded eigenvalues. The radiation condition is formulated as
$$
\lim_{R\to\infty}\frac{1}{\log R}\int_{1/R<y<1}\big\|\big(D_y - i\sigma_{\pm}(k)\big)u(y)\big\|_{L^2(M_0)}^2\frac{dy}{y^n} = 0,
$$
$$
\lim_{R\to\infty}\frac{1}{\log R}\int_{1<y<R}\big\|\big(D_y - i\sigma_{\mp}(k)\big)u(y)\big\|_{L^2(M_0)}^2\frac{dy}{y^n} = 0,
$$
$$
\sigma_{\pm}(k) = \frac{n-1}{2} \mp ik.
$$
Note that these are  standard radiation conditions

(\ref{S3Outgoingrad}) and (\ref{S3Incomingrad}) since
$r=|{\log y}|$.

The inverse scattering from regular ends works well as in \S 3.
Namely, suppose we are given two such asymptotically hyperbolic manifolds $\mathcal M^{(1)}$, $\mathcal M^{(2)}$ whose  regular ends $\mathcal M_1^{(1)}$ and $\mathcal M_1^{(2)}$ are isometric, and the $(1,1)$ components of the S-matrix coincide for all $k$. Then $\mathcal M^{(1)}$ and $\mathcal M^{(2)}$ are isometric.
This is proven in \cite{IK12}.

We remark  that S{\'a} Barreto \cite{SaBa05}, using the framework of scattering metric due to Melrose, proved the existence of isometry between two asymptotically hyperbolic manifolds having the same scattering matrix without the assumption that $\mathcal M_1^{(1)}$ and $\mathcal M_1^{(2)}$ are isometric, although the decay assumption at infinity are different from ours. In \cite{GuSaBa08}, it is extended to asymptotically hyperbolic complex manifolds. 

\subsection{Poisson integral and the space $\mathcal B^{\ast}$}
Before entering into the problem of cusp, we briefly look at some aspects of Poisson integral. Let us take the most basic example of Poincar{\'e} disc : $D = \{z \in {\mathbb C}\, ; \, |z|<1\}$. It is well-known that the Poisson integral
\begin{equation}
u(z) = \frac{1}{2\pi}\int_0^{2\pi}\left(\frac{1-|z|^2}{|e^{i\theta}-z|^2}\right)^s
f(\theta)d\theta,
\label{S4PoissonDisc}
\end{equation}
$f(\theta)$ being a function on the boundary $\partial D = S^1$, is a solution to the equation
\begin{equation}
(-\Delta_G - 4s(1-s))u = 0
\label{S4DiscHelmholtz}
\end{equation}
for all $s \in {\mathbb C}$. 
Our space $\mathcal B^{\ast}$, for which $f(\theta) \in L^2(S^1)$,  has the following special feature : $\{u \in \mathcal B^{\ast}\, ; \, (-\Delta_G-\frac{1}{4}-k^2)u=0\}$ is the smallest space of solutions as regard to decay at infinity. In fact, the Rellich type theorem says that the solution of (\ref{S4DiscHelmholtz}) decaying faster than the elements of $\mathcal B^{\ast}$ vanishes identically. Although this space is the smallest, it contains sufficiently many solutions for the inverse scattering. The largest solution space was found by Helgason \cite{Hel}, namely all solutions of (\ref{S4DiscHelmholtz}) are represented as (\ref{S4PoissonDisc}), where $f(\theta)$ is Sato's hyperfunction. This theorem was extended to the general symmetric spaces by Kashiwara-Kowata-Minemura-Okamoto-Oshima-Tanaka \cite{KKMOOT}.
The theorem of Helgason suggests that one may control the solution space of the equation $(-\Delta_G-\frac{1}{4}-k^2)u=0$ through function spaces on the boundary at infinity.
This enables us to extend the notion of S-matrix. We utilize this idea for the cusp.

\section{Arithmetic surface and generalized S-matrix}
We restrict ourselves here to the 2-dimensional case. Lots of examples of hyperbolic surfaces are given by the action of discrete groups on the upper half plane. If we consider the geometrically finite case, which implies that the associated quotient space  is a polygon with sides of finite number of geodesics, the infinity consists of a finite number of cusp and funnel (the latter being a slightly perturbed  regular infinity discussed in \ref{S4HypSpace}). However, if the Fuchsian group contains elliptic elements, the quotient space $\mathcal M$ has singularities. Note that $\mathcal M$ itself is an analytic manifold without singularities, however, the hyperbolic metric induced from ${\mathbb C}_+$ becomes singular at elliptic fixed points. 
Thus they are orbifolds, moreover, many examples of quotient spaces appearing in number theory have only cusp at infinity. They are non-compact, but has a finite volume. The result of the previous section cannot be applied to this case.
As a matter of fact,  the S-matrix at a cusp end does not have enough information to reconstruct the metric, since the cusp gives only one-dimensional contribution to the continuous spectrum.  A  remedy lies in generalizing the notion of S-matrix. 

Assume that $\mathcal M_1 = (1,\infty)\times (-1/2,1/2)$ with two sides identified.
Take any solution $u$ of the equation $(H-k^2)u=0$, and expand it into the Fourier series on $\mathcal M_1$:
$$
u(x,y) = \sum_{n\in{\bf Z}}e^{2\pi nx}u_n(y).
$$
Then $u_n$ satisfies
$$
y^2\big(-\partial_y^2 + (2ny)^2\big)u_n(y) - \frac{1}{4}u_n(y) =k^2u_n(y) ,
$$
hence is written as
$$
u_n(y) = 
\left\{
\begin{split}
&\tilde a_n\,y^{1/2} I_{-ik}(2\pi|n|y) + \tilde b_n\,y^{1/2}K_{ik}(2\pi|n|y), \quad n \neq 0, \\
&a_0 \, y^{1/2-ik} + b_0\, y^{1/2+ik}, \quad n=0.
\end{split}
\right.
$$
where $I_{\nu}$ and $K_{\nu}$ are modified Bessel functions. We then see that {\it all} solutions of the  equation $(H-k^2)u=0$ behave like
\begin{equation}
u(x,y) \sim a_0y^{1/2-ik} + \sum_{n\neq0}a_ne^{inx +|n|y} + b_0y^{1/2+ik} + \sum_{n\neq0}
b_ne^{inx-|n|y},
\label{S5Cuspexpgrowing}
\end{equation}
as $y \to \infty$. Given an exponentially growing wave 
represented by the first two terms in (\ref{S5Cuspexpgrowing}), one can uniquely construct 
a solution of $(H-k^2)u=0$ belonging to $\mathcal B^{\ast}$ on ends $\mathcal M_j$, $ j\neq 1$, 
and behaving like (\ref{S5Cuspexpgrowing}) in $\mathcal M_1$. We call the mapping
\begin{equation}
\mathcal S_{11}(k) : \{a_n\}_{n\in{\bf Z}} \to \{b_n\}_{n\in{\bf Z}}
\label{S5GeneralizedSmatrix}
\end{equation}
the {\it generalized S-matrix}. 
Passing to the Fourier series, we see that to define (\ref{S5GeneralizedSmatrix}), we are using a class of analytic functionals bigger than Sato's hyperfunction. Note that $\mathcal S_{11}(k)$ is an infinite matrix, and the usual S-matrix $S_{11}(k)$ is its $(0,0)$-entry. This generalized S-matrix has enough information for the inverse scattering. We make an artificial boundary in $\mathcal M_1$ and consider the boundary value problem in the interior domain. Then the knowledge of $\mathcal S_{11}(k)$ for all $k$ enables us to determine the N-D map $\Lambda_{int}(k)$, hence to apply the BC method.

In the interior domain $\mathcal M_{int}$, there is a finite number singular points. Since they have a special structure, we can deal with $\mathcal M_{int}$ as a manifold with {\it conical  singularities}, to which we can apply the BC method as well.
We can thus prove the following theorem.

\begin{theorem}\label{S5ThCusp}
Suppose we are given two 2-dimensional asymptotically hyperbolic surfaces with conical singularities $\mathcal M^{(1)}, \mathcal M^{(2)}$. Suppose they have pure cusp ends $\mathcal M_1^{(1)}, \mathcal M_1^{(2)}$, i.e. with metric $((dy)^2 + (dx)^2)/y^2$, and their generalized S-matrices coincide : $\mathcal S_{11}^{(1)}(k) = \mathcal S_{11}^{(2)}(k)$, for all $ k$. Then $\mathcal M^{(1)}$ and $\mathcal M^{(2)}$ are isometric.
\end{theorem}

 In particular, if $\mathcal H^{(i)} = \Gamma^{(i)}\slash {\bf H}^2$ with geometrically finite Fuchsian group $\Gamma^{(i)}$, and their generalized S-matrix coincide, $\Gamma^{(1)}$ and $\Gamma^{(2)}$ are conjugate. See \cite{IKL11} for details.
Also note that the results in subsections 4.5 and 5.1 also hold when $\mathcal M$ has both  
regular ends and cusps.

As in the previous cases, the key idea for the proof of Theorem \ref{S5ThCusp} is the use of exponentially growing solutions of the Helmholtz equation, which appears naturally here because of the form of hyperbolic metric. For cylindrical ends, we have encountered a similar situation, in which case, however, the non-physical scattering matrix is obtained by the analytic continuation of the physical scattering matrix. 
 This is not true for the cusp. In fact, Zelditch \cite{Ze} constructed an example of non-isometric hyperbolic surfaces with the same physical scattering matrix.

\section{Works in progress}
\subsection{Higher dimensional asymptotically hyperbolic orbifolds}
In higer dimensions, we have new phenomenon in the problem in \S 5. Consider the case of $n=3$, and let the Picard group $\Gamma = SL(2,{\bf Z}+i{\bf Z})$ act on ${\bf H}^3$ through quarternions. The resulting quotient space is a 3-dimensional hyperbolic manifold with singularities, i.e. an orbifold. Unlike the 2-dimensional case, the singularities are not confined in a compact set. They form an unbounded curve, and the manifold at infinity is not a smooth manifold, but a 2-dimensional orbifold (see \cite{ElGrMe}). However, the results for the forward problem also hold in this case with no essential change, and we are also expecting the inverse scattering results as well.

The structure of the manifold at infinity of general hyperbolic manifold is  complicated (see e.g. \cite{His94}). The forward and inverse scattering for this general case are challenging problems.

\subsection{Intermediate metrics}
We have worked on the problems (1)--(4), however, there is no reason that we must restrict ourselves to these cases. There is a wide area of problems for the intermediate metrics between (1) and (4), and those outside. For example, 
\cite{Kum10}, \cite{ItoSki} are dealing with non-existence of embedded eigenvalues, and \cite{Kum12} studies the LAP. We expect that the inverse scattering theory can also be developed in these cases.

\bibliographystyle{amsalpha}

\begin{thebibliography}{A}

\bibitem{Ag1} 
S. Agmon, \textit{Spectral theory of Schr{\"o}dinger operators on Euclidean and non-Euclidean spaces}, Comm. Pure Appl. Math. \textbf{39}, Supplement (1986), S 3 - S 16.

\bibitem{AKKLT} 
M. Anderson, A. Katsuda, Y. Kurylev, M. Lassas,  and M. Taylor, \textit{Boundary regularity for the Ricci equation, Geometric Convergence, and Gel'fand's Inverse Boundary Problem}, Inventiones Mathematicae {\bf 158} (2004), 261-321

\bibitem{Be87}
M. I. Belishev, \textit{An approach to multidimensional inverse problems for the wave equation}, Dokl. Akad. Nauk SSSR \textbf{297} (1987), 524-557.

\bibitem{Be97}
M. I. Belishev, \textit{Boundary control in construction of manifolds and metrics (the BC method)}, Inverse Problems \textbf{13} (1997), R1-R45. 

\bibitem{BeBla92}
M. I. Belishev and A. S. Blagovestcenskii, \textit{Multidimensional analogs of the Gel'fand-Levitan-Krein equations in inverse problems for the wave equations}, Ill-posed problems of Mathematical Physics and Analysis, Nobosibirsk : Nauk (1992), 50-63.

\bibitem{BeKu87}
M. I. Belishev, \textit{A nonstationary inverse problem for the multidimensional wave equation in the large}, ZAP. Nauchn. Sem. LOMI \textbf{165} (1987), 21-30.

\bibitem{BeKu92}
M. I. Belishev and Y. V. Kurylev, \textit{To the reconstruction of a Riemannian manifold via its spectral data (BC method)}, Comm. in P. D. E. \textbf{17} (1992), 767-804.

\bibitem{Bingham} K. Bingham, Y. Kurylev, M. Lassas and S. Siltanen,
\textit{Iterative time reversal control for inverse problems.} Inverse Problems and Imaging {\bf 2} (2008), 63-81.


\bibitem{Bla71a}
A. S. Blagovestcenskii, \textit{The local method of solution of the non-stationary inverse scattering problem for an inhomogeneous string}, Trudy Mat. Inst. Steklova \textbf{115} (1971), 28-38.


\bibitem{Bla71b}
A. S. Blagovestcenskii, \textit{The nonselfadjoint inverse matrix boundary problem for a hyperbolic differential equation}, in Problems of Mathematical Physics, \textbf{5},  Spectral Theory, Izdat. Leningrad Univ., Leningrad (1971), 38-62.

\bibitem{Borth}
D. Borthwick, \textit{Spectral Theory of Infinite-Area Hyperbolic Surfaces}, Progress in Mathematics \textbf{256}, Birkh{\"a}user, Boston-Basel-Berlin (2007).

\bibitem{BorJuPer05}
D. Borthwick, C. Judge and P. Perry, \textit{Selberg's zeta function and the spectral geometry of geometrically finite hyperbolic surfaces}, Comment. Math. Helv. \textbf{80} (2005), 483-515.


\bibitem{Tanya09}
T. Christiansen, \textit{Sojourn times, manifolds with infinite cylindrical ends, and inverse problem for planer waveguides}, J. Anal. Math. \textbf{107} (2009), 79-106.

\bibitem{Eid69}
D. M. Eidus, \textit{The principle of limit amplitude}, Russ. Math. Survey \textbf{24} (1969), 97-167.

\bibitem{ElGrMe}
J. Elstrodt, F. Grunewald and J. Mennicke, \textit{Groups Acting in Hyperbolic Spaces}, Springer, Berlin-Heidelberg-New York (1998).

\bibitem{Fa76}
L. D. Faddeev, \textit{Inverse problem of quantum scattering theory}, J. Sov. Math. \textbf{5} (1976), 334-396.




\bibitem{GeLe}
I. M. Gel'fand and B. M. Levitan, \textit{On the determination of a differential equation from its spectral function}, Izv. Akad. Nauk USSR, Ser. Mat. \textbf{15} (1951), 309-360.

\bibitem{GLU1}
 A. Greenleaf, M. Lassas, and G. Uhlmann, \textit{On nonuniqueness for Calderon's inverse problem}, Math. Res. Lett. {\bf 10} (2003),  685-693.

\bibitem{GKLU0}
A. Greenleaf, Y. Kurylev, M. Lassas, and G. Uhlmann, \textit{Full-wave invisibility of active devices at all frequencies.} Commun. in Math. Phys. {\bf 275} (2007), 749-789

\bibitem{GKLU1}
 A. Greenleaf, Y. Kurylev, M. Lassas, and G. Uhlmann,
  \textit{Invisibility and Inverse Problems}. Bulletin of the Amer. Math. Soc. {\bf 46} (2009), 55-97.

\bibitem{GKLU2}
 A. Greenleaf, Y. Kurylev, M. Lassas, and G. Uhlmann:,  \textit{Cloaking Devices, Electromagnetic Wormholes and Transformation Optics}. SIAM Review {\bf 51} (2009), 3-33. 

\bibitem{GuSaBa08}
C. Guillarmou and A. S{\'a} Barreto, \textit{Scattering and inverse scattering on ACH manifolds}, J. reine angew. Math. \textbf{622} (2008), 1-55.

\bibitem{Hel}
S. Helgason, \textit{A duality for symmetric spaces with applications to group representations}, Adv. in Math. \textbf{5} (1970), 1-154.

\bibitem{His94}
P. D. Hislop, \textit{The geometry and spectra of hyperbolic manifolds}, Proc. Indian Acad. Sci. (Math. Sci.) \textbf{104} (1994), 715-776.

\bibitem{IsaNa93}
V. Isakov and A. Nachman, \textit{Global uniqueness for a two-dimensional semilinear elliptic inverse problem}, Trans. Amer. Math. Soc., \textbf{347} (1993), 3375-3390.

\bibitem{Iso04}
H. Isozaki, \textit{Inverse spectral problems on hyperbolic manifolds and their applications to inverse boundary value problems in Euclidean space}, Amer. J. Math. \textbf{126} (2004), 1261-1313. 

\bibitem{IK12}
H. Isozaki and Y. Kurylev, \textit{Introduction to spectral theory and inverse problems on asymptotically hyperbolic manifolds}, arXiv 1102.5382.

\bibitem{IKL10}
H. Isozaki, Y. Kurylev and M. Lassas, \textit{Forward and inverse scattering problem on manifolds with asymptotically cylindrical ends}, Journal of Funct. Anal. \textbf{258} (2010), 2060-2118.

\bibitem{IKL11}
H. Isozaki, Y. Kurylev and M. Lassas, \textit{Conic singularities, generalized scattering matrix and inverse scattering on asymptotically hyperbolic surfaces}, arXiv:1108.1577.

\bibitem{Iso03}
H. Isozaki, \textit{Inverse spectral theory}, in \textit{Topics in the Theory of Schr{\"o}dinger Operators}, edit. H. Araki and H. Ezawa, World Scientific (2003), 93-143.

\bibitem{ItoSki}
K. Ito and E. Skibsted, \textit{Absence of embedded eigenvalues for Riemannian Laplacians}, arXiv 1109.1928.


\bibitem{KKL01}
A. Katchalov, Y. Kurylev and M. Lassas, \textit{Inverse Boundary Spectral Problems}, Chapman and Hall/CRC, Monographs and Surveys in Pure and Applied Mathematics, \textbf{123} (2001).

\bibitem{KKLM}
 A. Katchalov, Y. Kurylev, M. Lassas, and N. Mandache, \textit{Equivalence of time-domain inverse problems and boundary spectral problem}, Inverse problems {\bf 20} (2004), 419-436

\bibitem{KKMOOT}
M. Kashiwara, A. Kowata, K. Minemura, K. Okamoto, T. Oshima, M. Tanaka, \textit{Eigenfunctions of invariant differential operators on a symmetric spaces}, 
Ann. Math. \textbf{107} (1978), 1-39.

\bibitem{Ke} M.\ Kerker,  {\it Invisible bodies}, J. Opt. Soc. Am. {\bf 65}
(1975), 376.


\bibitem{KheNov}
G. M. Khenkin and R. G. Novikov, \textit{The $\overline{\partial}$ equation in the multi-dimensional inverse scattering problem}, Russian Math. Surveys, \textbf{42} (1987), 109-180.

\bibitem{Kr51a}
M. G. Krein, \textit{Determination of the density of an inhomogeneous string from its spectrum}, Dokl. Akad. Nauk SSSR \textbf{76} (1951), 345-348.

\bibitem{Kr51b}
M. G. Krein, \textit{On inverse problems for an inhomogeneous string}, Dokl. Akad. Nauk SSSR \textbf{82} (1951), 669-672.

\bibitem{KruL} K. Krupchyk, Y. Kurylev, M. Lassas, \textit{Inverse spectral problems on a closed manifold}.
{Journal de Mathematique Pures et Appliquees}
{\bf 90} (2008), 42-59. 
  

\bibitem{Kum10}
H. Kumura, \textit{The radial curvature of an end that makes eigenvalues vanish  in the essential spectrum, I}, Math. Ann. \textbf{346} (2010), 795-828.

\bibitem{Kum12}
H. Kumura, \textit{Limiting absorption principle on manifolds having ends with various measure growth rate limits}, arXiv. math DG/0606125.


\bibitem{KLY} Y. Kurylev, M. Lassas, T. Yamaguchi. \textit{Uniqueness and stability in inverse
spectral problems for collapsing manifolds,} arxiv.math: 1209.587v2


\bibitem{LSU}
 M. Lassas, M. Salo, G. Uhlmann, \textit{Wave imaging}. Chapter in Handbook of Mathematical Methods in Imaging, O. Scherzer (Ed.), Springer-Verlag, 2011, 50 pp

\bibitem{LTU}
 M. Lassas, M. Taylor and G. Uhlmann, \textit{The Dirichlet-to-Neumann map for complete Riemannian manifolds with boundary}, Comm. in Anal. and Geom. {\bf 11} (2003), 207-222.

\bibitem{La-U} M. Lassas and G. Uhlmann, \textit{ On determining a Riemannian 
manifold from the Dirichlet-to-Neumann map},  Annales Scientifiques 
de l' Ecole Normale Superi\'eure, {\bf 34} (2001), 771--787. 

\bibitem{Mar55}
V.  A. Marchenko, \textit{The construction of the potential energy from the phases of the scattering waves}, Dokl. Akad. Nauk SSSR \textbf{104} (1955), 695-698.

\bibitem{MaMe87}
R. Mazzeo and R. Melrose, \textit{Meromorphic extension of the resolvent on complete spaces with asymptotically constant negative curvature}, J. Funct. Anal. \textbf{75} (1987), 260-310.

\bibitem{MaVa05}
R. Mazzeo and A. Vasy, \textit{Analytic continuation of the resolvent of the Laplacian on symmetric spaces of noncompact type}, J. Funct. Anal. \textbf{228} (2005), 311-368.

\bibitem{Me95}
R. Melrose, \textit{Geometric Scattering Theory}, Cambridge University Press (1995).


\bibitem{MU}
R. Melrose, G. Uhlmann, 
{\textit Generalized Backscattering and the Lax-Phillips Transform}, 
Serdica Math. J. {\bf 34} (2008), 1026-1044.

\bibitem{Mou81}
E. Mourre, \textit{Absence of singular continuous spectrum for certain self-adjoint operators}, Comm. in Math. Phys. \textbf{78} (1981), 391-400.

\bibitem{Mu87}
W. M{\"u}ller, \textit{Manifolds with cusp of rank one, spectral theory and $L^2$-index theorem}, Lecture Notes in Math. \textbf{1224}, Springer (1987).

\bibitem{Na}
A. Nachman, \textit{Reconstruction from boundary measurements}, Ann. of Math. \textbf{143} (1988), 531-576.



\bibitem{PSU}
L. P\"aiv\"arinta, M. Salo, and G. Uhlmann, \textit{Inverse scattering for the magnetic Schr{\"o}dinger operator}, 
J. Funct. Anal. {\bf 259} (2010), 1771-1798.

\bibitem{Perry07}
P. A. Perry, \textit{The spectral geometry of geometrically finite hyperbolic manifolds}, Spectral Theory and Mathematical Physics: A Festschrift in Honor of Barry Simon's 60th Birthday, Proc. Sympos. Pure Math., \textbf{76}, A. M. S. Providence, RI (2007), 289-328.

\bibitem{ReSi}
M. Reed and B. Simon, \textit{Methods of Modern Mathematical Physics, III : Scattering Theory}, Academic Press (1979).

\bibitem{SaBa05}
A. S{\'a} Barreto, \textit{Radiation fields, scattering and inverse scattering on asymptotically hyperbolic manifolds}, Duke Math. J. \textbf{129} (2005), 407-480.

\bibitem{SylUhl87} 
J. Sylvester and G. Uhlmann, \textit{A global uniqueness theorem for an inverse boundary value problem}, Ann. of Math. \textit{125} (1987), 153-169.

\bibitem{Tata}
D. Tataru, \textit{Unique continuation for solutions to PDEs; between H{\"o}rmander's theorem and Holmgren's theorem}, 
Comm. in PDE. \textbf{20} (1995), 855-884.

\bibitem{Tata2}
D. Tataru, \textit{On the regularity of boundary traces for the wave equation.}
{Ann. Scuola Norm. Sup. Pisa Cl. Sci.}  {\bf 26}  (1998), 185--206.


\bibitem{Ya}
D. R. Yafaev, \textit{Mathematical Scattering Theory, Analytic Theory}, Mathematical Surveys and Monographs \textbf{158}, A. M. S., Providence, Rhode Island (2009).

\bibitem{Ze}
S. Zelditch, \textit{Kuznecov sum formulae and Szeg{\"o} limit formulae on manifolds}, Commun. in P. D. E. \textbf{17} (1992), 221-260.


\end{thebibliography}

\end{document}